\documentclass[12pt]{article}
\usepackage{theorem}
\usepackage{amsmath}
\usepackage{amssymb}
\usepackage{epsfig}
\usepackage{graphics}

\newsavebox{\coextremal}

\newtheorem{thm}{Theorem}[section]
\newtheorem{prop}[thm]{Proposition}

{\theorembodyfont{\rm}

\newtheorem{rem}[thm]{Remark}}
{\theorembodyfont{\bf}
}

\newcommand{\ra}{\rightarrow}

\newcommand{\Seg}{\mathrm{Seg}}
\newcommand{\NS}{\mathrm{NS}}
\newcommand{\NE}{\mathrm{NE}}
\newcommand{\NM}{\mathrm{NM}}
\newcommand{\N}{\mathrm{N}}
\newcommand{\Mov}{\mathrm{Mov}}

\newcommand{\PGL}{\mathrm{PGL}}

\newcommand{\bP}{{\mathbb P}}

\newcommand{\bQ}{{\mathbb Q}}
\newcommand{\bR}{{\mathbb R}}
\newcommand{\bZ}{{\mathbb Z}}

\newcommand{\cJ}{{\mathcal J}}

\newcommand{\cL}{{\mathcal L}}

\newcommand{\fS}{{\mathfrak S}}
\newcommand{\fh}{{\mathfrak h}}

\newcommand{\oM}{{\overline M}}

\title{On the effective cone of the moduli space of pointed
rational curves}
\author{Brendan Hassett and Yuri Tschinkel}
\date{November 22, 2001}
\begin{document}
\maketitle

\section{Introduction}

For a smooth projective variety, Kleiman's criterion for
ample divisors states that the closed ample cone (i.e., the nef cone)
is dual to the closed cone of effective curves.  Since the work of
Mori, it has been clear that {\em extremal rays} of the cone of
effective curves play a special role in birational geometry.  
These correspond to certain distinguished supporting hyperplanes of
the nef cone which are negative with respect to the canonical class.  
Contractions of extremal rays are the fundamental
operations of the minimal model program.    

Fujita \cite{F} has initiated a dual theory, with the (closed) cone of
effective divisors playing the central role.  It is natural then to
consider the dual cone and its generators.  Those which are
negative with respect to the canonical  class are called {\em coextremal rays},
and have been studied by Batyrev \cite{Ba}.  They are expected to play
a fundamental role in Fujita's program of classifying fiber-space structures
on polarized varieties.  

There are relatively few varieties for which the extremal and 
coextremal rays are fully understood.  Recently, moduli spaces of pointed
rational curves $\oM_{0,n}$ have attracted considerable
attention, especially in connection with mathematical physics and 
enumerative geometry.  
Keel and McKernan  first considered the
`Fulton conjecture':  The cone of effective curves of $\oM_{0,n}$
is generated by one-dimensional boundary strata.  This is
proved for $n\le 7$ \cite{KeMc}.  The analogous statement for divisors,
namely, that the effective cone of $\oM_{0,n}$ is generated
by boundary divisors, is known to be false (\cite{Ke} and \cite{Ve}).
The basic idea is to consider the map
$$
r:\oM_{0,2g}\hookrightarrow \oM_g, \quad n=2g,
$$
identifying pairs $(i_1i_2),(i_3i_4),\ldots,(i_{2g-1}i_{2g})$
of marked points to nodes.  There exist effective divisors
in $\oM_g$ restricting to effective
divisors not spanned by boundary divisors (see Remark~\ref{rem:hyp}).
However, it is true that for each $n$ the cones of $\fS_n$-invariant 
effective divisors are generated by boundary 
divisors \cite{KeMc}.  
 
\

In recent years it has become apparent that various arithmetic questions
about higher dimensional algebraic varieties defined over
number fields 
are also closely related to the cone of effective divisors. 
For example, given a  variety $X$ over a number field $F$, 
a line bundle $L$  in the interior of $\NE^1(X)$,
an open $U\subset X$ over which $L^N (N\gg 0)$
is globally generated, and a height $H_{\cL}$ associated to 
some adelic metrization $\cL$ of $L$, we can consider the
asymptotic behavior of the counting function
$$
N(U,\cL,B)=\# \{x\in U(F)\,|\,  H_{\cL}(x)\le B\}\quad B>0.
$$
There is a heuristic principle that,
after suitably restricting $U$,
$$
N(U,\cL,B)= c(\cL) B^{a(L)} \log(B)^{b(L)-1}(1+o(1)),
$$
as $B\ra \infty$  (see \cite{BT}). Here 
$$
a(L):=\inf\{ a\in \bR\,|\, aL+K_X\in \NE^1(X)\},
$$
$b(L)$ is the codimension of the face of $\NE^1(X)$ containing
$a(L)L+K_X$ (provided that $\NE^1(X)$ is 
locally polyhedral at this point),
and $c(\cL)>0$ is a constant depending on the chosen height
(see \cite{BM} and \cite{BT} for more details).
Notice that the explicit 
determination of the constant $c(\cL)$
also involves the knowledge of the effective cone.

Such asymptotic formulas can be proved for
smooth complete intersections 
in $\bP^n$ of small degree using the classical circle method
in analytic number theory 
and for varieties closely 
related to linear algebraic groups, 
like flag varieties, 
toric varieties etc., using adelic harmonic
analysis (\cite{BT} and references therein). 
No general techniques to treat arbitrary varieties
with many rational points are currently available. 
To our knowledge, the only other variety for which
such asymptotic is known to hold is the moduli space
$\oM_{0,5}$ (Del Pezzo surface of degree 5) in its
anticanonical embedding \cite{DB}.
Upper and lower bounds, with the expected $a(L)$ and $b(L)$,
are known (see \cite{VW}) for the Segre cubic threefold
$$
\Seg=\{(x_0,\ldots,x_5):
\sum_{j=0}^5 x_j^3=\sum_{j=0}^5 x_j=0\}.
$$ 
This admits an explicit resolution
by the moduli space $\oM_{0,6}$ (Remark~\ref{rem:segre});  
see \cite{Hu} for the relationship between the Segre cubic
and moduli spaces.

\

Our main result (Theorem~\ref{thm:main}) is a computation of
the effective cone of $\oM_{0,6}$.  Besides the boundary
divisors, the generators are the loci in $\oM_{0,6}$
fixed under 
$$\sigma=(i_1i_2)(i_3i_4)(i_5i_6)\in \fS_6, \quad 
\{i_1,i_2,i_3,i_4,i_5,i_6\}=\{1,2,3,4,5,6\}.$$
This equals the closure of $r^*\fh\cap M_{0,6}$,
where $\fh$ is the hyperelliptic locus in $\oM_3$.  
The effective and moving cones of $\oM_3$ are studied 
in detail by Rulla \cite{Ru}.  Rulla's inductive analysis
of the moving cone is similar to the method outlined 
in Section~\ref{sect:effgen}.  Results on the ample
cone of $\oM_{0,6}$ have been recently obtained by Farkas
and Gibney \cite{FG}.

\

The arithmetic consequences of Theorem~\ref{thm:main}
will be addressed in a future paper.

\

\noindent {\bf Acknowledgments:}  We are grateful to 
Se\'an Keel for helpful discussions, especially concerning
the moving cone of $\oM_{0,6}$, and to William Rulla, for
pointing out an error in an early version of this paper.           
The first author was partially supported by the 
Institute of Mathematical Sciences
of the Chinese University of Hong Kong and NSF grant 0196187.
The second author was partially supported by the NSA, NSF and the Clay 
Foundation.

\section{Generalities on effective cones}
\label{sect:effgen}

Let $X$ be a nonsingular projective variety with N\'eron-Severi
group $\NS(X)$ and group of one-cycles $\N_1(X)$.  
The closed effective cone of $X$ is the closed convex cone 
$$\NE^1(X)\subset \NS(X)\otimes \bR$$
generated by effective divisors on $X$.  Let $\NM_1(X)$ be
the dual cone $\NE^1(X)^*$ in $\N_1(X)\otimes \bR$.  Similarly,
let $\NE_1(X)$ be the cone of effective curves and $\NM^1(X)$ its
dual, the nef cone.  
  
We review one basic strategy for computing $\NE^1(X)$.  
Suppose we are given a collection 
$\Gamma=\{A_1,\ldots,A_m\}$
of effective divisors that we expect to generate the effective 
cone and a subset 
$
\Sigma\subset \Gamma.
$ 
For any effective divisor $E$, we have a decomposition
$$
E=M_{\Sigma}+B_{\Sigma}, \quad B_{\Sigma}
=a_1A_1+\ldots+a_mA_m,\quad a_j\ge 0,
$$ 
where $B_{\Sigma}$ is the fixed part of $|E|$ supported in $\Sigma$.  
The divisor $M_{\Sigma}$ may have fixed components, but they are
not contained in $\Sigma$.  If $\Mov(X)_{\Sigma}$ denotes the 
closed cone generated by effective divisors without fixed components 
in $\Sigma$, then $M_{\Sigma} \in \Mov(X)_{\Sigma}$.
Any divisor of $\Mov(X)_{\Sigma}$ restricts to an effective
divisor on each $A_j\in \Sigma$.  Consequently,
$$\Mov(X)_{\Sigma}\subset \NM_1(\Sigma,X)^*,$$
where $\NM_1(\Sigma,X)\subset \N_1(X)$ is generated by the images
of the $\NM_1(A_i)$.  To prove that $\Gamma$ generates
$\NE^1(X)$, it suffices then to check that
$$\{\text{cone generated by $\Gamma$}\}^*\subset \NM_1(\Sigma,X).$$

\section{Geometry of $\oM_{0,n}$}

\subsection{A concrete description of $\oM_{0,n}$}
\label{subsect:concrete}

In this section we give a basis for the N\'eron-Severi 
group of $\oM_{0,n}$ and write down the boundary divisors 
and the symmetric group action.

\

We recall the explicit iterated blow-up realization 
$$
\beta_n:\oM_{0,n} \ra {\bP}^{n-3}
$$
from \cite{Has} (see also a related construction in \cite{Kap}.)
This construction involves choosing one of the marked
points;  we choose $s_n$.  
Fix points $p_1,\ldots,p_{n-1}$ in linear general 
position in $\bP^{n-3}:=X_0[n]$.  
Let $X_1[n]$ be the blow-up of $\bP^{n-3}$ at $p_1,\ldots,p_{n-1}$,
and let $E_1,\ldots,E_{n-1}$ denote the exceptional divisors
(and their proper transforms in subsequent blow-ups).  
Consider the proper transforms $\ell_{ij}\subset X_1[n]$ of the lines
joining $p_i$ and $p_j$.  Let $X_2[n]$ be the blow-up of $X_1[n]$
along the $\ell_{ij}$, with exceptional divisors $E_{ij}$.
In general, $X_k[n]$ is obtained from $X_{k-1}[n]$
by blowing-up along proper transforms of the linear spaces spanned by
$k$-tuples of the points.  The exceptional divisors are denoted
$$E_{i_1,\ldots,i_k}\quad \{i_1,\ldots,i_k\} \subset \{1,\ldots,n-1\}.$$
This process terminates with a nonsingular variety $X_{n-4}[n]$
and a map 
$$
\beta_n\,:\,X_{n-4}[n]\ra \bP^{n-3}.
$$  
One can prove that
$X_{n-4}[n]$ is isomorphic to $\oM_{0,n}$.     
We remark that for a generic point $p_n\in \bP^{n-3}$, we have
an identification
$$
\beta_n^{-1}(p_n)=(C,p_1,p_2,\ldots,p_n),
$$
where $C$ is the unique rational normal curve of degree $n-3$
containing $p_1,\ldots,p_n$ (see \cite{Kap} for further
information).

Let $L$ be the pull-back
of the hyperplane class on ${\bP}^{n-3}$ by $\beta_n$.  We obtain
the following explicit basis for $\NS(\oM_{0,n})$:
$$
\{ L ,E_{i_1},E_{i_1i_2}, \ldots ,E_{i_1,\ldots,i_k},\ldots,E_{i_1,\ldots, i_{n-4}}\}.
$$
We shall use the following dual basis for the one-cycles 
$\N_1(\oM_{0,n})$:
$$
\{ L^{n-4},(-E_{i_1})^{n-4},\ldots,
(-E_{i_1,\ldots,i_k})^{n-3-k}L^{k-1},
\ldots,(-E_{i_1,\ldots,i_{n-4}})L^{n-5}\}. \quad (\dagger)
$$

\subsection{Boundary divisors}

Our next task is to identify the boundary divisors of $\oM_{0,n}$
in this basis.  These are indexed by partitions
$$\{1,2,\ldots,n\}=S\cup S^c, \quad n\in S \text{ and }|S|,|S^c|\ge 2;$$
the generic point of the divisor $D_S$ corresponds to a curve
consisting of two copies of ${\bP}^1$ intersecting at a node $\nu$, 
with marked points from $S$ on one component and from $S^c$
on the other.  Thus we have an isomorphism 
\begin{eqnarray*}
D_S &\simeq& \oM_{0,|S|+1}\times \oM_{0,|S^c|+1}, \hskip 2cm (\ddag) \\
(\bP^1,S)\cup_{\nu}(\bP^1,S^c) & \longrightarrow & 
(\bP^1,S\cup\{\nu \})\times (\bP^1,S^c\cup\{\nu \}).
\end{eqnarray*}

The exceptional divisors are identified as follows:
$$E_{i_1,\ldots,i_k}=D_{i_1,\ldots,i_k,n},  \quad \{i_1,\ldots,i_k\}
\subset \{1,\ldots,n-1\}, k\le n-4.$$
The remaining divisors $D_{i_1,\ldots,i_{n-3},n}$ are the proper transforms of the 
hyperplanes spanned by $(n-3)$-tuples of points;  we have
$$[D_{i_1,\ldots,i_{n-3},n}]=L-E_{i_1}-E_{i_2}-\ldots-
E_{i_1,\ldots,i_{n-4}}-\ldots-E_{i_2,\ldots,i_{n-3}}.$$     

\begin{rem}\label{rem:segre}
The explicit resolution of the Segre threefold
$$R:\oM_{0,6}\ra \Seg$$
alluded to in the introduction
is given by the linear series 
$$
|2L-E_1-E_2-E_3-E_4-E_5|.
$$
The image is a cubic threefold with ten ordinary double
points, corresponding to the lines $\ell_{ij}$ contracted
by $R$. 
\end{rem}

\subsection{The symmetric group action on $\oM_{0,n}$}
The symmetric group $\fS_n$ acts
on $\oM_{0,n}$ by the rule
$$\sigma(C,s_1,\ldots,s_n)=(C,s_{\sigma(1)},\ldots,s_{\sigma(n)}).$$
Let $F_{\sigma}\subset M_{0,n}$ denote the locus 
fixed by an element $\sigma\in \fS_n$.

We make explicit the $\fS_n$-action in terms of our blow-up realization.
Choose coordinates $(z_0,z_1,z_2,\ldots,z_{n-3})$ on $\bP^{n-3}$ so that 
$$
p_1=(1,0,\ldots,0),  \ldots,  p_{n-2}=(0,\ldots,0,1), 
\ p_{n-1}=(1,1,\ldots,1,1).
$$
Each permutation of the first $(n-1)$ points can be realized by a unique 
element of $\PGL_{n-2}$.  For elements of $\fS_n$ fixing $n$,
the action on $\oM_{0,n}$ is induced by the corresponding linear transformation
on $\bP^{n-3}$.  Now let $\sigma=(jn)$ and consider the commutative
diagram
$$\begin{array}{ccc}
\oM_{0,n} & \stackrel{\sigma}{\ra} &  \oM_{0,n} \\
\beta_n \downarrow \quad &  	& \quad \downarrow \beta_n \\
\bP^{n-3} &   \stackrel{\sigma'}{\dashrightarrow} & \bP^{n-3}  
\end{array}.$$
The birational map $\sigma'$ is the Cremona transformation based at the 
points $p_{i_1},\ldots,p_{i_{n-2}}$ where 
$$\{i_1,\ldots,i_{n-2},j\}=\{1,2,\ldots,n-1\},$$ 
e.g., when $\sigma=(n-1,n)$ we have
$$\sigma(z_0,z_1,\ldots,z_{n-3})=
(z_1z_2\ldots z_{n-3},z_0z_2\ldots z_{n-3},\ldots,z_0\ldots z_{n-4}).$$

\section{Analysis of surfaces in $\oM_{0,6}$}

\subsection{The $\oM_{0,5}$ case}
\begin{prop}\label{prop:five}
$\NE^1(\oM_{0,5})$ is generated by the divisors
$D_{ij},$ where $\{ij\}\subset \{1,2,3,4,5\}$.  
\end{prop}    

\noindent
{\em Sketch proof:}
This is well-known, but we sketch the basic ideas to introduce notation
we will require later.  As we saw in \S~\ref{subsect:concrete},
$\oM_{0,5}$ is the blow-up of $\bP^2$ at four points in general position.
Consider the set of boundary divisors
$$\Sigma=\{D_{i5},D_{ij}\}=\{E_i,L-E_i-E_j\},
\quad \{i,j\} \subset \{1,2,3,4\}$$
and the set of semiample divisors
$$\Xi=\{L-E_i,2L-E_1-E_2-E_3-E_4,L,2L-E_i-E_j-E_k\},
\ \{i,j,k\}\subset \{1,2,3,4\}.$$
These semiample divisors come from the forgetting maps
$$\phi_i:\oM_{0,5} \ra \oM_{0,4}\simeq \bP^1, \quad i=1,\ldots,5$$
and the blow-downs
$$\beta_i:\oM_{0,5}\ra \bP^2,\quad i=1\ldots,5.$$
Kleiman's criterion yields
$$C(\Sigma) \subset \NE_1(\oM_{0,5})=\NM^1(\oM_{0,5})^*\subset C(\Xi)^*.$$
All the inclusions are equalities because
the cones generated by $\Xi$ and $\Sigma$ are dual;  this
can be verified by direct computation (e.g., using 
the computer program {\tt PORTA}~\cite{PORTA}).  
$\square$

\subsection{Fixed points and the Cayley cubic}
\label{subsect:Cayley}

We identify the fixed-point {\em divisors} for 
the $\fS_6$-action on $\oM_{0,6}$.
When $\tau=(12)(34)(56)$ we have
$$\tau(z_0,z_1,z_2,z_3)=(z_0z_2z_3,z_1z_2z_3,z_0z_1z_2,z_0z_1z_3)$$
and $F_{\tau}$ is given by $z_0z_1=z_2z_3.$
It follows that
$$[F_{\tau}]=2L-E_1-E_2-E_3-E_4-E_5-E_{13}-E_{23}-E_{24}-E_{14}.$$   
More generally, when $\tau=(ab)(cd)(j6)$ we have
$$[F_{\tau}]=2L-E_1-E_2-E_3-E_4-E_5-E_{ac}-E_{ad}-E_{bc}-E_{bd}.$$
\begin{rem}\label{rem:hyp}
Consider $(\bP^1,s_1,\ldots,s_6)\in F_{\tau}$ and
the quotient under the corresponding involution
$$q:\bP^1\longrightarrow \bP^1,\quad q(s_1)=q(s_2),\
q(s_3)=q(s_4),\text{ etc.}$$
Consider the map
$r:\oM_{0,6} \ra \oM_3$
identifying the pairs $(12),(34),$ and $(56)$ and
write $C=q(\bP^1,s_1,\ldots,s_6)$, so there is an induced
$q':C \ra \bP^1.$
Thus $C$ is hyperelliptic
and $F_{\tau}$ corresponds to the closure of $r^*\fh\cap\oM_{0,6}$,
where $\fh\subset \oM_3$ is the hyperelliptic locus.  
\begin{figure}[hhh]
\centerline{\psfig{figure=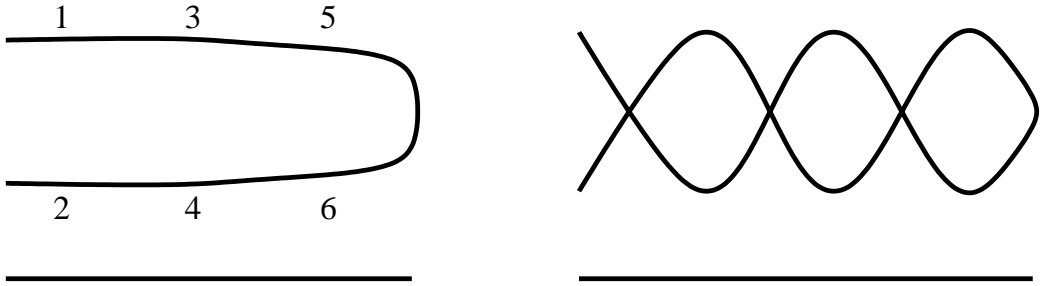}}
\caption{Trinodal hyperelliptic curves}
\end{figure}
\end{rem}

It will be useful to know the effective cone of the
fixed point divisors $F_{\sigma}$.  We have seen that these
are isomorphic to ${\bP}^1\times {\bP}^1$ blown-up at five
points $p_1,\ldots,p_5$.  The projection from $p_5$
$$\bP^3 \dashrightarrow \bP^2$$
induces a map $\varphi:F_{\sigma}\ra \bP^2$,
realizing $F_{\sigma}$ as a blow-up of $\bP^2$:  
Take four general lines $\ell_1,\ldots,\ell_4$  
in $\bP^2$ with intersections $q_{ij}=\ell_i\cup \ell_j$, and
blow-up $\bP^2$ along the $q_{ij}$.  We write
$$
\NS(F_{\sigma})={\bZ}H+{\bZ}G_{12}+\ldots+{\bZ}G_{34},
$$
where the $G_{ij}$ are the exceptional divisors and $H$ is
the pull back of the hyperplane class from $\bP^2$.    

\begin{prop}\label{prop:fixcone} 
$\NE^1(F_{\sigma})$ is generated by the $(-1)$-curves
$$G_{12},\ldots,G_{34},H-G_{ij}-G_{kl},$$
and the $(-2)$-curves
$$H-G_{ij}-G_{ik}-G_{il},\quad \{i,j,k,l\}=\{1,2,3,4\}.$$
\end{prop}

\noindent {\em Proof:}
Let $\Sigma$ be the above collection of 13 curves.  Consider also
the following collection $\Xi$ of $38$ divisors, grouped as 
orbits under the $\fS_4$-action:
$$ 
\begin{array}{c|c|l}
\text{typical member} & \text{orbit size} & \text{induced morphism} \\
\hline
H & 1 & \text{blow-down $\varphi:F_{\sigma}\ra \bP^2$} \\  
H-G_{12} & 6 & \text{conic bundle $F_{\sigma}\ra \bP^1$} \\
2H-G_{12}-G_{13}-G_{23} & 4 & \text{blow-down $F_{\sigma}\ra \bP^2$}\\
2H-G_{12}-G_{23}-G_{34} & 12 & \text{blow-down $F_{\sigma}\ra \bP^2$}\\
2H-G_{12}-G_{23}-G_{34}-G_{14} & 3 & \text{conic bundle $F_{\sigma}\ra \bP^1$}\\
3H-2G_{12}-G_{13}-G_{23}-G_{34} & 12 & \text{blow-down $F_{\sigma}\ra \bP(1,1,2)$}
\end{array}
$$
Note that each of these divisors is semiample:  the corresponding
morphism is indicated in the table.  In particular, 
\begin{eqnarray*}
C(\Sigma)&:=&\{\text{cone generated by $\Sigma$}\} \subset \NE_1(F_{\sigma}), \\
C(\Xi)&:=&\{\text{cone generated by $\Xi$}\} \subset \NM^1(F_{\sigma}) 
\end{eqnarray*}
and Kleiman's criterion yields
$$C(\Sigma) \subset \NE_1(F_{\sigma})=\NM^1(F_{\sigma})^*\subset C(\Xi)^*.$$
A direct verification using {\tt PORTA}~\cite{PORTA} shows that the cones
$C(\Sigma)$ and $C(\Xi)$ are dual, so all the inclusions are equalities.
$\square$

\begin{rem}
The image of $F_{\tau}$ under the resolution $R$
of \ref{rem:segre} is a cubic surface with four double points,
classically called the {\em Cayley cubic} \cite{Hu}.  
\end{rem}

\section{The effective cone of $\oM_{0,6}$}
We now state the main theorem:
\begin{thm}\label{thm:main}
The cone of effective divisors $\NE^1(\oM_{0,6})$  
is generated by the boundary divisors and the fixed-point divisors
$F_{\sigma}$, where $\sigma\in \fS_6$ is a product of
three disjoint transpositions.
\end{thm}

\subsection{Proof of Main Theorem}
We use the strategy outlined in \S~\ref{sect:effgen}.  
Consider the collection of boundary and fixed-point loci
$$\Gamma=\{D_{ij},D_{ijk},F_{\sigma}, \quad \sigma=(ij)(kl)(ab),
\quad \{i,j,k,l,a,b\}=\{1,2,3,4,5,6\}\}$$
and the subset of boundary divisors
$$\Sigma=\{D_{ij},D_{ijk}\}.$$
We take 
$$\Xi=\{ \text{ images $\rho \in \N_1(\oM_{0,6})$ of 
generators of $\NM_1(A_i), A_i\in \Sigma$ }\}.$$  

We compute the cone $\NM_1(\Sigma,\oM_{0,6})$,
the convex hull of the union of the images of
$\NM_1(D_{ij})$ and $\NM_1(D_{ijk})$ in $\N_1(\oM_{0,6})$.
Throughout, we use the dual basis for $\N_1(\oM_{0,6})$
(cf. $(\dag)$):
\begin{eqnarray*}
& &\{L^2,E_1^2,E_2^2,E_3^2,E_4^2,E_5^2,-LE_{12},-LE_{13},\\
& &-LE_{14},-LE_{15},-LE_{23},-LE_{24},-LE_{25},-LE_{34},
-LE_{35},-LE_{45}  \}.
\end{eqnarray*}
Recall the isomorphism (\ddag)
$$(\pi_{ijk},\pi_{lab}):D_{ijk} \longrightarrow 
\bP^1\times \bP^1, \quad \{i,j,k,l,a,b,c\}=\{1,2,3,4,5,6\} $$
so that
$$\N_1(D_{ijk})=\bZ B_{ijk} \oplus \bZ B_{lab}, \quad
\NM_1(D_{ijk})=\bR_+ B_{ijk}+ \bR_+B_{lab},$$
where $B_{ijk}$ is the class of the fiber of $\pi_{ijk}$.  
For example, the inclusion $j_{345}:D_{345} \hookrightarrow \oM_6$
induces
$$
(j_{345})_*=
\left( \begin{array}{cccccccccccccccc}
0 & 0 & 0 & 0 & 0 & 0 & -1 & 0 & 0 & 0 & 0 & 0 & 0 & 0 & 0 & 0 \\
1 & 1 & 1 & 0 & 0 & 0 & -1 & 0 & 0 & 0 & 0 & 0 & 0 & 0 & 0 & 0
\end{array} \right)^T
$$
using the bases $(\dagger)$ for ${\rm N}_1(\oM_{0,6})$ and 
$\{B_{126},B_{345}\}$ for $D_{345}$.  In particular, we find
$$\NM_1(\{D_{ijk}\},\oM_{0,6}) = C(\{ B_{ijk} \}),
\quad \{i,j,k\}\subset \{1,2,3,4,5,6\},$$
with $\binom{6}{3}=20$ generators permuted transitively
by $\fS_6$ (Table~\ref{table:Bijk}).
 
\begin{table}[p]
\caption{Generators for $\NM_1(\{D_{ijk}\},\oM_{0,6})$}
\label{table:Bijk}
\scriptsize
$$
\begin{array}{c|cccccccccccccccc}
\hline
B_{126}&0&0&0&0&0&0&-1&0&0&0&0&0&0&0&0&0 \\
B_{136}&0&0&0&0&0&0&0&-1&0&0&0&0&0&0&0&0 \\
B_{146}&0&0&0&0&0&0&0&0&-1&0&0&0&0&0&0&0 \\
B_{156}&0&0&0&0&0&0&0&0&0&-1&0&0&0&0&0&0 \\
B_{236}&0&0&0&0&0&0&0&0&0&0&-1&0&0&0&0&0 \\
B_{246}&0&0&0&0&0&0&0&0&0&0&0&-1&0&0&0&0 \\
B_{256}&0&0&0&0&0&0&0&0&0&0&0&0&-1&0&0&0 \\
B_{346}&0&0&0&0&0&0&0&0&0&0&0&0&0&-1&0&0 \\
B_{356}&0&0&0&0&0&0&0&0&0&0&0&0&0&0&-1&0 \\
B_{456}&0&0&0&0&0&0&0&0&0&0&0&0&0&0&0&-1 \\
B_{123}&1&0&0&0&1&1&0&0&0&0&0&0&0&0&0&-1 \\
B_{124}&1&0&0&1&0&1&0&0&0&0&0&0&0&0&-1&0 \\
B_{125}&1&0&0&1&1&0&0&0&0&0&0&0&0&-1&0&0 \\
B_{134}&1&0&1&0&0&1&0&0&0&0&0&0&-1&0&0&0 \\
B_{135}&1&0&1&0&1&0&0&0&0&0&0&-1&0&0&0&0 \\
B_{145}&1&0&1&1&0&0&0&0&0&0&-1&0&0&0&0&0 \\
B_{234}&1&1&0&0&0&1&0&0&0&-1&0&0&0&0&0&0 \\
B_{235}&1&1&0&0&1&0&0&0&-1&0&0&0&0&0&0&0 \\
B_{245}&1&1&0&1&0&0&0&-1&0&0&0&0&0&0&0&0 \\
B_{345}&1&1&1&0&0&0&-1&0&0&0&0&0&0&0&0&0
\end{array}
$$
\normalsize
\end{table}

The boundary divisor $D_{ij}$ is isomorphic
to $\oM_{0,5}$ with marked points $\{k,l,a,b,\nu\}$    
where $\{i,j,k,l,a,b\}=\{1,2,3,4,5,6\}$
and $\nu$ is the node (cf. formula (\ddag)).    
By Proposition~\ref{prop:five}, the cone $\NM_1(D_{ij},\oM_{0,6})$ 
is generated by the classes
$$\{ A_{ij},A_{ij;k},A_{ij;l},A_{ij;a},A_{ij;b},
C_{ij},C_{ij;k},C_{ij;l},C_{ij;a},C_{ij;b} \}\subset \N_1(\oM_{0,6})$$
corresponding to the forgetting and blow-down morphisms
$$\{
\phi_{\nu},\phi_k,\phi_l,\phi_a,\phi_b, 
\beta_{\nu},\beta_k,\beta_l,\beta_a,\beta_b 
\}.$$

As an example, consider the inclusion
$j_{45}:D_{45}\hookrightarrow \oM_{0,6}$ with
$${j_{45}}_*=\left(\begin{array}{cccccccccccccccc}
1&0&0&0&0&0&  1&1&0&0&1&0&0&0&0&0\\
0&1&0&0&0&0& -1&-1&0&0&0&0&0&0&0&0\\
0&0&1&0&0&0& -1&0&0&0&-1&0&0&0&0&0\\
0&0&0&1&0&0&  0&-1&0&0&-1&0&0&0&0&0\\
0&0&0&0&0&0&  0&0&0&0&0&0&0&0&0&1
\end{array} \right)^T.
$$
Applying Prop.~\ref{prop:five},
we obtain generators for $\NM_1(D_{45},\oM_{0,6})$
(Table~\ref{table:D45}).
\begin{table}
\caption{Generators for $\NM_1(D_{45},\oM_{0,6})$} 
\label{table:D45}
\scriptsize
$$
\begin{array}{c|cccccccccccccccc}
\hline
A_{45}& 	1&0&0&0&0&0 &1&1&0&0&1&0&0&0&0&1 \\
A_{45;1}& 	1&1&0&0&0&0 &0&0&0&0&1&0&0&0&0&0 \\ 
A_{45;2}&	1&0&1&0&0&0 &0&1&0&0&0&0&0&0&0&0 \\
A_{45;3}&	1&0&0&1&0&0 &1&0&0&0&0&0&0&0&0&0 \\
A_{45;6}&	2&1&1&1&0&0 &0&0&0&0&0&0&0&0&0&1 \\
C_{45;6}&	1&0&0&0&0&0 &1&1&0&0&1&0&0&0&0&0 \\
C_{45;1}&	2&0&1&1&0&0 &1&1&0&0&0&0&0&0&0&1 \\
C_{45;2}&	2&1&0&1&0&0 &1&0&0&0&1&0&0&0&0&1 \\
C_{45;3}&	2&1&1&0&0&0 &0&1&0&0&1&0&0&0&0&1 \\
C_{45}&	2&1&1&1&0&0 &0&0&0&0&0&0&0&0&0&0 \\
\end{array}
$$
\normalsize
\end{table}

Quite generally, four $(-1)$-curves in $D_{ij}$ are contained
in $D_{ijk},D_{ijl},D_{ija},$ and $D_{ijb}$,
with classes $B_{ijk},B_{ijl},B_{ija},$ and $B_{ijb}$ respectively.    
Thus we have the relations
$$
C_{ij}=A_{ij;k}+B_{ijk}, \quad C_{ij;k}=A_{ij}+B_{ijk}
$$
which implies that the $C_{ij}$ and $C_{ij;k}$ are redundant:

\begin{prop}
The cone $\NM_1(\Sigma,\oM_{0,6})$ is generated by the $A_{ij}$, the $A_{ij;k}$,
and the $B_{ijk}$.  
\end{prop}    
These are written out in Tables \ref{table:Bijk},\ref{table:Aij}, and \ref{table:Aijk}.
\begin{table}
\caption{Generators $A_{ij}$ for $\NM_1(\{D_{ij}\},\oM_{0,6})$} 
\label{table:Aij}
\scriptsize
$$
\begin{array}{c|cccccccccccccccc}
\hline
A_{12}& 1&  0&0&0&0&0&  1&0&0&0&0&0&0&1&1&1 \\
A_{13}& 1&  0&0&0&0&0&  0&1&0&0&0&1&1&0&0&1 \\
A_{14}& 1&  0&0&0&0&0&  0&0&1&0&1&0&1&0&1&0 \\
A_{15}& 1&  0&0&0&0&0&  0&0&0&1&1&1&0&1&0&0 \\
A_{16}& 0&  -2&0&0&0&0& 1&1&1&1&0&0&0&0&0&0 \\
A_{23}& 1&  0&0&0&0&0&  0&0&1&1&1&0&0&0&0&1 \\
A_{24}& 1&  0&0&0&0&0&  0&1&0&1&0&1&0&0&1&0 \\
A_{25}& 1&  0&0&0&0&0&  0&1&1&0&0&0&1&1&0&0 \\
A_{26}& 0&  0&-2&0&0&0& 1&0&0&0&1&1&1&0&0&0 \\
A_{34}& 1&  0&0&0&0&0&  1&0&0&1&0&0&1&1&0&0 \\
A_{35}& 1&  0&0&0&0&0&  1&0&1&0&0&1&0&0&1&0 \\
A_{36}& 0&  0&0&-2&0&0& 0&1&0&0&1&0&0&1&1&0 \\
A_{45}& 1&  0&0&0&0&0&  1&1&0&0&1&0&0&0&0&1 \\
A_{46}& 0&  0&0&0&-2&0& 0&0&1&0&0&1&0&1&0&1 \\
A_{56}& 0&  0&0&0&0&-2& 0&0&0&1&0&0&1&0&1&1
\end{array}
$$
\normalsize
\end{table}

\begin{table}
\vspace{-.6in}

\caption{Generators $A_{ij;k}$ for $\NM_1(\{D_{ij}\},\oM_{0,6})$} 
\label{table:Aijk}

\scriptsize
$$\begin{array}{c|cccccccccccccccc}
\hline
A_{12;3}&1&0&0&1&0&0&   0&0&0&0&0&0&0&0&0&1\\
A_{12;4}&1&0&0&0&1&0&   0&0&0&0&0&0&0&0&1&0\\
A_{12;5}&1&0&0&0&0&1&   0&0&0&0&0&0&0&1&0&0\\
A_{12;6}&2&0&0&1&1&1&   1&0&0&0&0&0&0&0&0&0\\
A_{13;2}&1&0&1&0&0&0&   0&0&0&0&0&0&0&0&0&1\\
A_{13;4}&1&0&0&0&1&0&   0&0&0&0&0&0&1&0&0&0\\
A_{13;5}&1&0&0&0&0&1&   0&0&0&0&0&1&0&0&0&0\\
A_{13;6}&2&0&1&0&1&1&   0&1&0&0&0&0&0&0&0&0\\
A_{14;2}&1&0&1&0&0&0&   0&0&0&0&0&0&0&0&1&0\\
A_{14;3}&1&0&0&1&0&0&   0&0&0&0&0&0&1&0&0&0\\
A_{14;5}&1&0&0&0&0&1&   0&0&0&0&1&0&0&0&0&0\\
A_{14;6}&2&0&1&1&0&1&   0&0&1&0&0&0&0&0&0&0\\
A_{15;2}&1&0&1&0&0&0&   0&0&0&0&0&0&0&1&0&0\\
A_{15;3}&1&0&0&1&0&0&   0&0&0&0&0&1&0&0&0&0\\
A_{15;4}&1&0&0&0&1&0&   0&0&0&0&1&0&0&0&0&0\\
A_{15;6}&2&0&1&1&1&0&   0&0&0&1&0&0&0&0&0&0\\
A_{16;2}&0&-1&0&0&0&0&  1&0&0&0&0&0&0&0&0&0\\
A_{16;3}&0&-1&0&0&0&0&  0&1&0&0&0&0&0&0&0&0\\
A_{16;4}&0&-1&0&0&0&0&  0&0&1&0&0&0&0&0&0&0\\
A_{16;5}&0&-1&0&0&0&0&  0&0&0&1&0&0&0&0&0&0\\
A_{23;1}&1&1&0&0&0&0&   0&0&0&0&0&0&0&0&0&1\\
A_{23;4}&1&0&0&0&1&0&   0&0&0&1&0&0&0&0&0&0\\
A_{23;5}&1&0&0&0&0&1&   0&0&1&0&0&0&0&0&0&0\\
A_{23;6}&2&1&0&0&1&1&   0&0&0&0&1&0&0&0&0&0\\
A_{24;1}&1&1&0&0&0&0&   0&0&0&0&0&0&0&0&1&0\\
A_{24;3}&1&0&0&1&0&0&   0&0&0&1&0&0&0&0&0&0\\
A_{24;5}&1&0&0&0&0&1&   0&1&0&0&0&0&0&0&0&0\\
A_{24;6}&2&1&0&1&0&1&   0&0&0&0&0&1&0&0&0&0\\
A_{25;1}&1&1&0&0&0&0&   0&0&0&0&0&0&0&1&0&0\\
A_{25;3}&1&0&0&1&0&0&   0&0&1&0&0&0&0&0&0&0\\
A_{25;4}&1&0&0&0&1&0&   0&1&0&0&0&0&0&0&0&0\\
A_{25;6}&2&1&0&1&1&0&   0&0&0&0&0&0&1&0&0&0\\
A_{26;1}&0&0&-1&0&0&0&  1&0&0&0&0&0&0&0&0&0\\
A_{26;3}&0&0&-1&0&0&0&  0&0&0&0&1&0&0&0&0&0\\
A_{26;4}&0&0&-1&0&0&0&  0&0&0&0&0&1&0&0&0&0\\
A_{26;5}&0&0&-1&0&0&0&  0&0&0&0&0&0&1&0&0&0\\
A_{34;1}&1&1&0&0&0&0&   0&0&0&0&0&0&1&0&0&0\\
A_{34;2}&1&0&1&0&0&0&   0&0&0&1&0&0&0&0&0&0\\
A_{34;5}&1&0&0&0&0&1&   1&0&0&0&0&0&0&0&0&0\\
A_{34;6}&2&1&1&0&0&1&   0&0&0&0&0&0&0&1&0&0\\
A_{35;1}&1&1&0&0&0&0&   0&0&0&0&0&1&0&0&0&0\\
A_{35;2}&1&0&1&0&0&0&   0&0&1&0&0&0&0&0&0&0\\
A_{35;4}&1&0&0&0&1&0&   1&0&0&0&0&0&0&0&0&0\\
A_{35;6}&2&1&1&0&1&0&   0&0&0&0&0&0&0&0&1&0\\
A_{36;1}&0&0&0&-1&0&0&  0&1&0&0&0&0&0&0&0&0\\
A_{36;2}&0&0&0&-1&0&0&  0&0&0&0&1&0&0&0&0&0\\
A_{36;4}&0&0&0&-1&0&0&  0&0&0&0&0&0&0&1&0&0\\
A_{36;5}&0&0&0&-1&0&0&  0&0&0&0&0&0&0&0&1&0\\
A_{45;1}&1&1&0&0&0&0&   0&0&0&0&1&0&0&0&0&0\\
A_{45;2}&1&0&1&0&0&0&   0&1&0&0&0&0&0&0&0&0\\
A_{45;3}&1&0&0&1&0&0&   1&0&0&0&0&0&0&0&0&0\\
A_{45;6}&2&1&1&1&0&0&   0&0&0&0&0&0&0&0&0&1\\
A_{46;1}&0&0&0&0&-1&0&  0&0&1&0&0&0&0&0&0&0\\
A_{46;2}&0&0&0&0&-1&0&  0&0&0&0&0&1&0&0&0&0\\
A_{46;3}&0&0&0&0&-1&0&  0&0&0&0&0&0&0&1&0&0\\
A_{46;5}&0&0&0&0&-1&0&  0&0&0&0&0&0&0&0&0&1\\
A_{56;1}&0&0&0&0&0&-1&  0&0&0&1&0&0&0&0&0&0\\
A_{56;2}&0&0&0&0&0&-1&  0&0&0&0&0&0&1&0&0&0\\
A_{56;3}&0&0&0&0&0&-1&  0&0&0&0&0&0&0&0&1&0\\
A_{56;4}&0&0&0&0&0&-1&  0&0&0&0&0&0&0&0&0&1
\end{array}$$
\normalsize
\end{table}

Our next task is to write out the generators for the dual
cone $C(\Gamma)^*$, as computed by {\tt PORTA} \cite{PORTA}.  
Since $\Gamma$ is stable under the $\fS_6$ action, 
so are $C(\Gamma)$ and its dual cone.
For the sake of brevity, we only write $\fS_6$-representatives
of the generators, ordered by anticanonical degree.  

\begin{table}
\caption{$\fS_6$-orbits of coextremal rays of $\oM_{0,6}$} 
\label{table:coextremal}
$$\begin{array}{rc|c|cccccccccccccccc}
&\deg_{-K} & \text{order} & & & & & & & & & & & & & & & & \\
\hline
  & & & & & & & & & & & & & & & & & & \\
(1)&2&1&
3&0&0&0&0&0&1&1&1&1&1&1&1&1&1&1\\
(2)&2&6&
1&0&0&0&0&1&0&0&0&0&0&0&0&0&0&0\\
(3)&2&15&
2&0&0&0&0&0&0&0&1&1&0&1&1&1&1&0\\
(4)&2&45&
1&0&0&0&0&0&0&0&0&0&0&0&1&1&0&0\\
(5)&3&60&
1&0&0&0&0&0&0&0&0&0&0&0&0&0&0&1\\
(6)&3&72&
2&0&0&0&0&0&0&0&1&1&1&0&1&1&0&0\\
(7)&3&120&
2&0&0&0&0&1&0&0&0&1&0&1&0&1&0&0\\
(8)&3&120&
2&0&0&0&0&1&0&0&1&0&1&0&0&1&0&0\\
(9)&3&180&
2&0&0&0&0&0&0&0&0&1&0&1&1&1&1&0\\
(10)&4&6&
1&0&0&0&0&0&0&0&0&0&0&0&0&0&0&0\\
(11)&4&10&
3&0&0&1&1&1&2&0&0&0&0&0&0&0&0&0\\
(12)&4&30&
2&0&0&0&0&0&0&0&0&1&0&0&1&0&1&1\\
(13)&4&60&
3&0&0&0&0&0&0&0&0&2&1&1&1&1&1&1\\
(14)&4&90&
3&0&0&0&0&1&0&0&1&1&0&1&1&2&0&0\\
(15)&4&90&
3&0&0&0&0&2&0&1&1&0&1&1&0&0&0&0\\
(16)&4&180&
2&0&0&0&0&1&0&0&0&0&0&1&0&1&0&0\\
(17)&4&180&
3&0&0&0&0&0&0&0&1&1&1&0&2&2&0&1\\
(18)&4&360&
2&0&0&0&0&0&0&0&0&1&0&0&1&1&0&1\\
(19)&4&360&
3&0&0&0&0&1&0&0&1&1&1&0&1&2&0&0\\
(20)&4&360&
3&0&0&0&0&1&0&1&1&0&1&2&0&0&1&0\\
(21)&5&120&
2&0&0&0&0&0&0&0&0&0&0&0&1&0&1&1\\
(22)&5&360&
3&0&0&0&0&0&0&0&0&2&0&1&1&2&0&1\\
(23)&5&360&
4&0&0&0&0&2&0&1&1&1&2&2&0&0&0&0\\
(24)&6&360&
4&0&0&0&0&0&0&0&0&3&0&2&1&2&1&1\\
(25)&6&360&
5&0&0&0&0&2&0&1&2&1&2&3&0&0&1&0\\
 & &    & & & & & & & & & & & & & & & \\
\hline
 & &3905& & & & & & & & & & & & & & &  
\end{array}$$
\end{table}

The discussion of Section~\ref{sect:effgen} shows
that Theorem~\ref{thm:main} will follow from the inclusion
$$
C(\Gamma)^*\subset \NM_1(\Sigma,\oM_{0,6}).
$$
We express each generator of $C(\Gamma)^*$
as a sum (with positive coefficients) 
of the $\{A_{ij},A_{ij;k},B_{ijk}\}$.  
Both cones are stable under the $\fS_6$-action, 
so it suffices to produce expressions for one
representative of each $\fS_6$-orbit.  We use
the representatives from Table~\ref{table:coextremal}:
\begin{eqnarray*}
(1)&=&A_{15}+A_{13}+A_{35}+2B_{246} \quad
(2)=A_{34;5}+B_{126} \\
(3)&=&A_{15}+A_{14}+2B_{236} \quad
(4)=A_{25}+B_{146}+B_{136}\\
(5)&=&A_{23}+B_{146}+B_{156}+B_{236}\quad
(6)=A_{15}+A_{14}+B_{236}+B_{246}+B_{356} \\
(7)&=&A_{13;5}+A_{15}+B_{236}+B_{246}\quad
(8)=A_{12;5}+A_{14}+B_{256}+B_{356}\\
(9)&=&A_{24}+A_{34}+B_{126}+B_{136}+B_{156}\\
(10)&=&A_{25}+B_{136}+B_{146}+B_{256}+B_{346}\\
(11)&=&A_{34;5}+A_{35;4}+A_{25;3}+B_{146}\quad
(12)=A_{12}+A_{34}+2B_{126}+2B_{346}\\
(13)&=&A_{15}+A_{14}+A_{23}+2B_{146}+2B_{236}\\
(14)&=&A_{23;5}+A_{15}+A_{25}+B_{136}+B_{146}+B_{236}\\
(15)&=&A_{23;5}+A_{24;5}+A_{15}+B_{156}+B_{346}\\
(16)&=&A_{24;5}+A_{15}+B_{136}+B_{156}+B_{236}\\
(17)&=&A_{23}+2A_{25}+2B_{136}+2B_{146}\\
(18)&=&A_{12;3}+A_{34}+B_{126}+B_{136}+A_{36;1}\\
(19)&=&A_{12;5}+A_{15}+A_{25}+B_{136}+B_{246}+B_{346}\\
(20)&=&A_{13;5}+A_{35}+A_{45}+2B_{126}+B_{456}\\
(21)&=&A_{12}+A_{13}+B_{126}+B_{136}+B_{246}+B_{346}+B_{456}\\
(22)&=&A_{15}+A_{23}+A_{34}+B_{126}+B_{146}+B_{156}+2B_{236}\\
(23)&=&A_{13;5}+A_{14;5}+A_{23}+A_{13}+B_{256}+2B_{456}\\
(24)&=&A_{15}+A_{23}+A_{24}+A_{34}+B_{126}+B_{136}+B_{146}+B_{156}+2B_{236}\\
(25)&=&2A_{14}+A_{24}+2A_{13;5}+2B_{256}+2B_{356}\\
\end{eqnarray*}
This completes the proof of Theorem~\ref{thm:main}.$\square$
 
\subsection{Geometric interpretations of coextremal rays}
By definition, a {\em coextremal ray} 
$\bR_+\rho \subset \NM_1(X)$ satisfies the following
\begin{itemize}
\item{for any nontrivial $\rho_1,\rho_2 \in \NM_1(X)$
with $\rho_1+\rho_2\in \bR_+\rho$, 
$\rho_1,\rho_2\in \bR_+\rho$;}
\item{$K_X\rho<0$.}
\end{itemize}
Batyrev (\cite{Ba}, Theorem 3.3) shows that, for smooth
(or $\bQ$-factorial terminal) threefolds, the 
minimal model program yields a geometric interpretation
of coextremal rays.  They arise from diagrams
$$
\begin{array}{ccl}
X & \stackrel{\psi}{\dashrightarrow} & Y \\
  & 		    & \downarrow \mu \\
  &		    & B
\end{array}
$$
where $\psi$ is a sequence of birational contractions and
$\mu$ is a Mori fiber space.  The coextremal ray
$\rho=\psi^*[C]$, where
$C$ is a curve lying in the general fiber of $\mu$.  
These interpretations will hold for higher-dimensional varieties,
provided the standard conjectures of the minimal model program
are true.  

It is natural then to write down these Mori fiber space structures 
explicitly.  Our analysis makes reference to the list of
orbits of coextremal rays in Table~\ref{table:coextremal}:
\begin{enumerate}
\item[(1)]{The anticanonical series $|-K_{\oM_{0,6}}|$ yields
a birational morphism
$$\oM_{0,6} \ra \cJ\subset \bP^4$$
onto a singular quartic hypersurface, called the
{\em Igusa quartic} \cite{Hu}.  The conics $C\subset \cJ$
pull back to the coextremal ray.}
\item[(2)]{Forgetting any of the six marked points
$$\oM_{0,6} \ra \oM_{0,5}$$
yields a Mori fiber space, and the fibers are
coextremal.}
\item[(3)]{We define a conic bundle structure on $\oM_{0,6}$
by explicit linear series, using the blow-up description
of Subsection~\ref{subsect:concrete}.  Consider the 
cubic hypersurfaces in $\bP^3$ passing through the points and lines
$$p_1,p_2,p_3,p_4,p_5,\ell_{14},\ell_{15},\ell_{24},\ell_{25},
\ell_{34},\ell_{35}.$$
We can compute the projective dimension
$$\dim|3L-E_1-E_2-E_3-E_4-E_5-E_{14}-E_{15}-E_{24}-E_{25}-E_{34}-E_{35}|=2.$$
This series yields a conic bundle structure
$$\mu:\oM_{0,6} \ra \bP^2$$
collapsing the two-parameter family of conics passing through
the six lines above.}
\item[(4)]{For any two disjoint subsets 
$\{i,j\},\{k,l\}\subset \{1,2,3,4,5,6\}$
we consider the forgetting maps
$$\phi_{ij}:\oM_{0,6} \ra \bP^1, \quad \phi_{kl}:\oM_{0,6} \ra \bP^1.$$
Together, these induce a conic bundle structure
$$(\phi_{ij},\phi_{kl}):\oM_{0,6} \ra \bP^1 \times \bP^1.$$
The class of a generic fiber is coextremal.}
\end{enumerate}

\subsection{The moving cone}
Our analysis gives, implicitly, the moving cone of $\oM_{0,6}$:

\begin{thm}
The closed moving cone of $\oM_{0,6}$ is equal to 
$\NM_1(\Gamma,\oM_{0,6})^*$, where $\Gamma$ is the set
of generators for $\NE^1(\oM_{0,6})$.  
\end{thm}

In the terminology of \cite{Ru}, 
the `inductive moving cone' equals the `moving cone'.    
We computed the ample cones to the boundaries $D_{ij}$ and $D_{ijk}$
and the fixed-point divisors $F_{\sigma}$ (Proposition~\ref{prop:fixcone}); 
this determines the moving cone completely.  
However, finding explicit {\em generators}
for the moving cone is a formidable computational problem.

\

\noindent
{\em Proof:}  
Recall that $\oM_{0,6}$ is a log Fano threefold:  
$-(K_{\oM_{0,6}}+\epsilon \sum_{ij}D_{ij})$ is ample for small $\epsilon>0$
\cite{KeMc}.  Using Corollary 2.16 of \cite{KeHu}, it follows that
$\oM_{0,6}$ is a `Mori Dream Space'.   The argument of Theorem 3.4.4 
of \cite{Ru} shows that an effective divisor on $\oM_{0,6}$
that restricts to an effective divisor on each generator $A_i\in \Gamma$
is in the moving cone.  $\square$

\begin{rem}
Our proof of Theorem~\ref{thm:main} uses the cone 
$\NM_1(\Sigma,\oM_{0,6})^*$, rather than the (strictly) smaller 
moving cone.  Of course, if the coextremal rays are in
$\NM_1(\Sigma,\oM_{0,6})$, {\em a fortiori} they are in
$\NM_1(\Gamma,\oM_{0,6})$.
\end{rem}

{}

\end{document}